%Paper: q-alg/9503009
%From: cquesne@ulb.ac.be (Quesne Christiane)
%Date: Fri, 17 Mar 1995 08:27:06 +0100 (MET)
%%%%%%%%%%%%%%%%%%%%%%%%%%%%%%%%%%%%%%%%%%%%%%%%%%%%%%%%%%%%%%%%%%%%%%%%%%%%
%%%%%%%%%%%%%%%%%%          PLAIN TEX          %%%%%%%%%%%%%%%%%%%%%%%%%%%%%
%%%%%%%%%%%%%%%%%%%%%%%%%%%%%%%%%%%%%%%%%%%%%%%%%%%%%%%%%%%%%%%%%%%%%%%%%%%%
\input preprint.sty

\font\tenof=msym10
\def\R{\hbox{\tenof R}}
\def\N{\hbox{\tenof N}}

\def\cqfd{\qquad\qquad\vrule height 4pt depth 2pt width 5pt}
\def\bar{\hbox{\rm ---}}
%
%=========================================================================
%
\title{On some nonlinear extensions of the angular momentum algebra}

\author{C Quesne\footnote\dag{Directeur de recherches
FNRS}\footnote\ddag {E-mail: cquesne@ulb.ac.be}}

\address{Physique Nucl\'eaire Th\'eorique et Physique Math\'ematique,
  Universit\'e Libre de Bruxelles, Campus de la Plaine CP229, Bd. du
  Triomphe, B1050 Brussels, Belgium}

\vfill
\shorttitle{Nonlinear extensions of angular momentum algebra}


%
%=========================================================================
%
\beginabstract
Deformations of the Lie algebras so(4), so(3,1), and e(3) that leave their
so(3) subalgebra undeformed and preserve their coset structure are
considered. It is shown that such deformed algebras are associative for any
choice of the deformation parameters. Their Casimir operators are obtained
and some of their unitary irreducible representations are constructed. For
vanishing deformation, the latter go over into those of the corresponding
Lie algebras that contain each of the so(3) unitary irreducible
representations at most once. It is also proved that similar deformations of
the Lie algebras su(3), sl(3,\R), and of the semidirect sum of an abelian
algebra t(5) and so(3) do not lead to associative algebras.
\endabstract
%
%===========================================================================
%
\section{Introduction}
In recent years, many works have been devoted to the study of deformations
and extensions of Lie algebras and their applications in various branches of
physics. Some of them are carried out in the mathematically well-defined
framework of quasitriangular Hopf algebras and deal with the so-called
quantum groups and $q$-algebras (Drinfeld~1986, Jimbo~1985). Others put less
emphasis on the coalgebra structure, which is often dropped completely, but
instead insist on preserving some other property of the Lie algebra that is
deformed. In this second category, one finds for instance some deformed
algebras that can be realized in terms of deformed creation and annihilation
operators (Fairlie and Zachos~1991, Fairlie and Nuyts~1994).\par
%
%---------------------------------------------------------------------------
%
In the same class, there are also deformed algebras that have a coset
structure $g_d = h + v_d$, and can be viewed as nonlinear extensions of an
ordinary Lie algebra $h$ (Ro\v cek~1991). This means that their generators
can be separated into the generators $E_i$ of $h$ and some operators
$E_{\alpha}$ transforming as a representation of $h$, and commuting among
themselves to give a function of the $E_i$'s only. In other words, they
satisfy the commutation relations
$$[E_i, E_j] = c_{ij}^k E_k \qquad [E_i, E_{\alpha}] = (\tau_i)_{\alpha}
       ^{\beta} E_{\beta} \qquad [E_{\alpha}, E_{\beta}] =
       f_{\alpha\beta}(E_i) \eqno(1.1)$$
where $c_{ij}^k$ and $(\tau_i)_{\alpha}^{\beta}$ are the structure constants
and some matrix representation of $h$ respectively, while
$f_{\alpha\beta}(E_i)$ are formal power series in the $E_i$'s. The latter
are constrained by the associativity requirement, i.e., Jacobi identities,
and by the condition that for some limiting values of the parameters, $g_d$
goes over into some Lie algebra $g$ with coset structure $g = h + v$. The
simplest example, corresponding to $h = \hbox{\rm u(1)}$, and $g =
\hbox{\rm su(2)}$ or su(1,1), has been discussed in detail
(Polychronakos~1990, Ro\v cek~1991).\par
%
%---------------------------------------------------------------------------
%
The interest of such constructions in the infinite-dimensional case was
already noted some years ago in the context of quantum field theory and
statistical physics models, where they are known as W-algebras
(Zamolodchikov~1986, Schoutens \etal 1989). More recently, finite versions
of these W-algebras were introduced by considering symplectic reductions of
finite-dimensional simple Lie algebras (Tjin~1992, de Boer and Tjin~1993).
It was shown in particular that the finite W$^{(2)}_3$-algebra, known as
$\overline{\hbox{\rm W}}^{(2)}_3$, is related to the above-mentioned simplest
example of deformed algebra $g_d$.\par
%
%--------------------------------------------------------------------------
%
The latter has also made recently its appearance in various physical
problems. Let us mention three of them. First, the deformed su(2) (or
su(1,1)) algebra may be considered as a dynamical symmetry algebra in some
quantum many-body models with symmetry-preserving Hamiltonians, such as those
occurring in quantum optics (Karassiov~1994, Karassiov and Klimov~1994 and
references quoted therein). Next, it is related to generalized
deformed parafermions, and through the introduction of a Fermi-like
oscillator Hamiltonian, provides a new algebraic description of the bound
state spectra of the Morse and P\H oschl-Teller potentials
(Quesne~1994). Moreover, by superposing these generalized deformed
parafermions with ordinary bosons, one gets some deformations of
parasupersymmetric quantum mechanics with new and nontrivial properties
(Beckers \etal 1995). Further, the algebra $\overline{\hbox{\rm W}}^{(2)}_3$
may be considered as the symmetry algebra of the two-dimensional anisotropic
harmonic oscillator with frequency ration 2:1 (Bonatsos \etal 1994).\par
%
%---------------------------------------------------------------------------
%
Motivated by these applications, we shall consider in the present paper
another class of examples of deformed algebras $g_d$, which should be
physically relevant. It corresponds to the case where the undeformed
subalgebra~$h$ is the angular momentum algebra so(3), and the deformed
subspace $v_d$ is spanned by the $2\lambda+1$ components of an so(3)
irreducible tensor of rank~$\lambda$. Special emphasis will be laid on the
vector ($\lambda=1$) and quadrupole ($\lambda=2$) cases, corresponding to
deformations of so(4) and su(3) (or of their noncompact or nonsemisimple
variants), respectively.\par
%
%---------------------------------------------------------------------------
%
In the following section, the relevant associativity conditions are
established. The vector and quadrupole cases are then studied in detail in
sections~3 and~4, respectively. Finally, section~5 contains the
conclusion.\par
%
%=========================================================================
%
\section{Nonlinear extensions of so(3)}
Let the generators $E_i$ and $E_{\alpha}$, introduced in the previous
section, be the spherical components $L_m^{\vphantom{\dagger}} = (-1)^m
L_{-m}^{\dagger}$, $m=+1,0,-1$, and $T^{\lambda}_{\mu} = (-1)^{\mu}
T^{\lambda\dagger}_{-\mu}$, $\mu = \lambda$, $\lambda-1$,
$\ldots$,~$-\lambda$, of an angular momentum operator and of an irreducible
tensor of integer rank~$\lambda$, respectively. In such a case, it is
advantageous to write eq.~(1.1) in a coupled commutator form as follows:
$$\eqalignno{[L, L]^1_m &= - \sqrt{2} L_m &(2.1a) \cr
    \bigl[L, T^{\lambda}\bigr]^{\Lambda}_M &= - \sqrt{\lambda (\lambda+1)}\,
       \delta_{\Lambda,\lambda}\, \delta_{M,\mu}\, T^{\lambda}_{\mu} &(2.1b)
       \cr
    \bigl[T^{\lambda}, T^{\lambda}\bigr]^{\Lambda}_M &= f^{\Lambda}_M(L)
       &(2.1c) \cr}$$
where $f^{\Lambda}(L)$ is an irreducible tensor of rank~$\Lambda$, whose
components can be written as formal power series in the vector operator~$L$.
The definition of coupled commutators and some of their properties are
reviewed in appendix~1. From eq.~(A1.2), it results that the values of
$\Lambda$ in eq.~(2.1c) are restricted to odd integers, $\Lambda=1$, 3,
$\ldots$,~$2\lambda-1$.\par
%
%--------------------------------------------------------------------------
%
It has been shown by Gaskell \etal (1978) that the number of linearly
independent irreducible tensors of rank~$\Lambda$, whose components are
monomials of degree~$n$ in~$L$, is equal to one if $n = \Lambda + 2k$, $k=0$,
1, 2,~$\ldots$, and zero otherwise. Hence the explicit form of the functions
$f^{\Lambda}_M(L)$ is given by
$$f^{\Lambda}_M(L) = \gamma_{\Lambda} g_{\Lambda}(\bi{L}^2) \Bigl[\cdots
    \bigl[[L\times L]^2 \times L\bigr]^3 \times\cdots \Bigr]^{\Lambda}_M
    \qquad \Lambda = 1,3,\ldots,2\lambda-1 \eqno(2.2)$$
where $\gamma_1$ is some real normalization constant, $\gamma_{\Lambda}=1$
for $\Lambda\ne1$,
$$g_{\Lambda}(\bi{L}^2) = \sum_{k=0}^{\infty} a^{(\Lambda)}_k \bi{L}^{2k}
    \qquad a^{(\Lambda)}_k \in \R \qquad a^{(1)}_0 = +1,0,\hbox{\rm or} -1
    \eqno(2.3)$$
is a formal power series in the scalar operator $\bi{L}^2 = \sum_m (-1)^m L_m
L_{-m}$, and the last factor on the right-hand side of~(2.2) is a ``stretched''
product of $\Lambda$ operators~$L$. The deformed algebra~$g_d$ is therefore
a $(2K+1)$-th degree algebra, where $K=0$, 1, 2,~$\ldots$, provided
$a^{(\Lambda)}_k = 0$ if $2k + \Lambda > 2K + 1$, and at least one
$a^{(\Lambda)}_k$ with $2k + \Lambda = 2K + 1$ is different from zero. As
$K=0$ corresponds to an ordinary Lie algebra, we shall henceforth refer to
$K$ as the deformation order.\par
%
%------------------------------------------------------------------------
%
We shall be concerned here with the cases where $T^{\lambda}$ is an
irreducible tensor of rank~1 (vector operator) or rank~2 (quadrupole
operator). In the former case, we shall denote $T^1_m$ by $A_m$, $m=+1$,
0,~$-1$; in eq.~(2.1c), $\Lambda$ then takes the single value $\Lambda=1$.
In the latter case, we shall denote $T^2_{\mu}$ by $Q_{\mu}$, $\mu=+2$,
$+1$, 0, $-1$,~$-2$; in eq.~(2.1c), $\Lambda$ then runs over the values
$\Lambda=1$ and $\Lambda=3$. For future use, it is helpful to point out some
relations among linearly dependent irreducible tensors:
$$\eqalignno{[L \times L]^0_0 &= - \case{1}{\sqrt{3}} \bi{L}^2 &(2.4a) \cr
    [L \times L]^1_m &= - \case{1}{\sqrt{2}} L_m &(2.4b) \cr
    \bigl[[L \times L]^2 \times L\bigr]^1_m &= - \case{2}{\sqrt{15}} \bi{L}^2
       L_m + \case{1}{2} \sqrt{\case{3}{5}} L_m &(2.4c) \cr
    \bigl[[L \times L]^2 \times L\bigr]^2_{\mu} &= - \sqrt{\case{3}{2}}\, [L
       \times L]^2_{\mu}. &(2.4d) \cr}$$\par
%
%------------------------------------------------------------------------
%
The deformed algebra generated by $L_m$ and $T^{\lambda}_{\mu}$ will be
associative provided commutators~(2.1) satisfy Jacobi identity. For three
irreducible tensors $T^{\lambda_1}$, $U^{\lambda_2}$, $V^{\lambda_3}$, of
ranks $\lambda_1$, $\lambda_2$,~$\lambda_3$, respectively, the latter can be
written in a coupled form, as shown in eq.~(A1.5). If $\bigl(T^{\lambda_1},
U^{\lambda_2}, V^{\lambda_3}\bigr) = (L, L, L)$, $\bigl(L, L,
T^{\lambda}\bigr)$, or $\bigl(L, T^{\lambda}, T^{\lambda}\bigr)$,
equation~(A1.5) is automatically satisfied, whereas if $\bigl(T^{\lambda_1},
U^{\lambda_2}, V^{\lambda_3}\bigr) = \bigl(T^{\lambda},  T^{\lambda},
T^{\lambda}\bigr)$, it leads to the set of conditions
$$\eqalignno{\sum_{\Lambda_{12}=1,3}^{2\lambda-1} &\left\{\delta_{
       \Lambda_{12},\Lambda_{23}} - 2 (-1)^{\lambda-\Lambda} U(\lambda
       \lambda \Lambda \lambda; \Lambda_{12} \Lambda_{23})\right\}
       \Bigl[T^{\lambda}, \bigl[T^{\lambda}, T^{\lambda}\bigr]^{\Lambda_{12}}
       \Bigr]^{\Lambda}_M = 0 \cr
    &\Lambda_{23} = 1, 3, \ldots, 2\lambda-1 \qquad \Lambda = \left|\lambda
       -\Lambda_{23}\right|, \left|\lambda - \Lambda_{23}\right| + 1, \ldots,
       \lambda + \Lambda_{23} &(2.5) \cr}$$
where $U(\lambda \lambda \Lambda \lambda; \Lambda_{12} \Lambda_{23})$ denotes
a Racah coefficient in unitary form (Rose~1957). Note that for simplicity's
sake, from now on we shall drop the component label of all irreducible
tensors.\par
%
%--------------------------------------------------------------------------
%
By using the numerical values of Racah coefficients, one finds that the set
of conditions~(2.5) reduces to a single independent condition
$$\bigl[A,[A,A]^1\bigr]^0 = 0 \eqno(2.6)$$
in the vector case, and to two independent conditions
$$\eqalignno{2\sqrt{2} \bigl[Q,[Q,Q]^1\bigr]^1 + \sqrt{7}
       \bigl[Q,[Q,Q]^3\bigr]^1 &=0 &(2.7a) \cr
    \bigl[Q,[Q,Q]^1\bigr]^3 - 2 \bigl[Q,[Q,Q]^3\bigr]^3 &=0 &(2.7b) \cr}$$
in the quadrupole one. In the next two sections, we shall determine whether
these identities are satisfied when the inner commutators are given by
eqs.~(2.1c), (2.2), and~(2.3).\par
%
%==========================================================================
%
\section{The vector case}
In the vector case, the deformed algebra~$g_d$ is defined by the commutation
relations
$$\eqalignno{[L,L]^1 &= - \sqrt{2} L \qquad [L,A]^{\Lambda} = - \sqrt{2}\,
       \delta_{\Lambda,1} A &(3.1a,b) \cr
    [A,A]^1 &= f^1(L) = - \sqrt{2}\, g(\bi{L}^2) L = - \sqrt{2}
       \biggl(\sum_{k=0}^{\infty} a_k \bi{L}^{2k} \biggr) L &(3.1c) \cr}$$
where the normalization constant $\gamma_1$ has been set equal to
$-\sqrt{2}$. A $K$-th order deformation corresponds to $a_K\ne0$, and
$a_{K+1} = a_{K+2} = \cdots = 0$.\par
%
%xxxxxxxxxxxxxxxxxxxxxxxxxxxxxxxxxxxxxxxxxxxxxxxxxxxxxxxxxxxxxxxxxxxxxxxxxxx
%
\subsection{Associativity condition}
The algebra $g_d$ is associative provided $A$ satisfies eq.~(2.6). By using
eqs.~(3.1) and~(A1.4), the latter can be rewritten as
$$\Bigl[\bigl[A, g(\bi{L}^2)\bigr]^1 \times L\Bigr]^0 = 0. \eqno(3.2)$$\par
%
%---------------------------------------------------------------------------
%
In the undeformed case where $g(\bi{L}^2) = a_0$, this condition is trivially
fulfilled. Then $g_d$ reduces to an ordinary Lie algebra~$g$. According to
whether $a_0 = +1$, $-1$, or~0, $g$ is the orthogonal algebra so(4), the
pseudo-orthogonal algebra so(3,1), or the Euclidian algebra e(3), which is a
semidirect sum of an abelian algebra t(3) and so(3) (Biedenharn~1961,
Naimark~1964, B\H ohm~1979).\par
%
%---------------------------------------------------------------------------
%
For a first-order deformation, condition~(3.2) reduces to
$$\Bigl[\bigl[A, \bi{L}^2\bigr]^1 \times L\Bigr]^0 = 0. \eqno(3.3)$$
{}From eqs.~(2.4a), (A1.4), (3.1b), and~(A1.1), one obtains
$$\bigl[A,\bi{L}^2\bigr]^1 = 2 \Bigl\{A + \sqrt{2}\, [L \times A]^1\Bigr\}.
    \eqno(3.4)$$
Moreover, standard tensor algebra and recoupling techniques (Rose~1957) show
that
$$\bigl[[L \times A]^1 \times L\bigr]^0 = \bigl[L \times [L \times A]^1
    \bigr]^0 = \bigl[[L \times L]^1 \times A\bigr]^0 = - \case{1}{\sqrt{2}}
    [L \times A|^0 \eqno(3.5)$$
where in the last step, use has been made of eq.~(2.4b). By introducing
eqs.~(3.4) and~(3.5) into the left-hand side of eq.~(3.3), one finds that the
latter is identically satisfied.\par
%
%-------------------------------------------------------------------------
%
It is now an easy task to show that if condition~(3.2) is fulfilled for a
$K$th-order deformation, then it is also satisfied for a ($K+1$)th-order one.
{}From eq.~(A1.4), one indeed obtains
$$\eqalignno{\Bigl[\bigl[A, \bi{L}^{2k+2}\bigr]^1 \times L\Bigr]^0 &=
       \left[\left[\bi{L}^{2k} \times \left[A,\bi{L}^2\right]^1\right]^1
       \times L\right]^0 + \biggl[\Bigl[\bigl[A,\bi{L}^{2k}\bigr]^1 \times
       \bi{L}^2\Bigr]^1 \times L\biggr]^0 \cr
    &= \bi{L}^{2k} \left[\left[A, \bi{L}^2\right]^1 \times L\right]^0 +
       \Bigl[\bigl[A, \bi{L}^{2k}\bigr]^1 \times L\Bigr]^0 \bi{L}^2. &(3.6)
       \cr}$$
Hence, if the relation
$$\Bigl[\bigl[A, \bi{L}^{2k}\bigr]^1 \times L\Bigr]^0 = 0 \eqno(3.7)$$
is identically satisfied for $k=K$, then the same is true for $k=K+1$. This
completes the proof by induction of the following result:\par
\bigskip
\noindent{\it Proposition 1.} For any choice of the deformation parameters
$a_k$, $k=1$, 2,~$\ldots$, equation~(3.1) defines an associative algebra
$g_d$, which is a deformed so(4), so(3,1), or~e(3) algebra according to
whether $a_0=+1$, $-1$, or~0.\par
\bigskip
\noindent{\it Remark.} The first-order deformation of so(4) has already
been encountered elsewhere. It is indeed the dynamical symmetry algebra of
a particle moving in a three-dimensional space with constant curvature
under the influence of a Coulomb potential (Higgs~1979, Leemon~1979,
Granovskii \etal 1992, de Vos and van Driel~1993). In such a case, the
deformation parameter~$a_1$ is related to the space curvature.\par
%
%xxxxxxxxxxxxxxxxxxxxxxxxxxxxxxxxxxxxxxxxxxxxxxxxxxxxxxxxxxxxxxxxxxxxxxxxxxx
%
\subsection{Casimir operators}
It is well known that the Lie algebras $g = \hbox{\rm so(4)}$, so(3,1),
and~e(3) have two independent Casimir operators, which may be written as
$$C_1 = a_0 \bi{L}^2 + \bi{A}^2 \qquad C_2 = \bi{L} \cdot \bi{A} \eqno(3.8)$$
where, as usual, $\bi{L} \cdot \bi{A}$ denotes the scalar product $\sum_m
(-1)^m L_m A_{-m}$. The purpose of the present subsection is to show that the
operators~(3.8) can be deformed so as to provide Casimir operators of the
deformed algebra $g_d$.\par
%
%---------------------------------------------------------------------------
%
The case of the second Casimir operator is easily solved. One finds the
following result:\par
\bigskip
\noindent{\it Proposition 2.} When going from $g = \hbox{\rm so(4)}$,
so(3,1), or~e(3) to $g_d$, defined in~(3.1), the operator $\bi{L} \cdot
\bi{A}$ remains a Casimir operator, which we shall denote by $C_{2d}$.\par
\bigskip
\noindent{\it Proof.} By using (A1.4), (3.1), and (A1.2), one obtains
$$\left[A, \bi{L} \cdot \bi{A}\right]^1 = - \case{1}{\sqrt{2}} [A, A]^1 +
    \case{1}{\sqrt{2}} f(\bi{L}^2) [L, L]^1 = 0. \eqno(3.9)$$
This completes the proof as, by construction, $\bi{L} \cdot \bi{A}$ commutes
with $L$.\cqfd\par
\bigskip
%
%--------------------------------------------------------------------------
%
The case of the first Casimir operator is more involved and actually has not
been solved in full generality. We conjecture that (i) for any choice of the
real constants~$a_k$ in eq.~(3.1), it is possible to find some real constants
$b_k$, $k=1$, 2,~$\ldots$, such that
$$C_{1d} = h(\bi{L}^2) + \bi{A}^2 \qquad \hbox{\rm where} \qquad h(\bi{L}^2)
    = \sum_{k=1}^{\infty} b_k \bi{L}^{2k} \eqno(3.10)$$
is a Casimir operator of the deformed algebra $g_d$, and (ii) for a
$K$th-order deformation, these constants are such that $b_k=0$ for $k>K+1$.
We shall now proceed to prove that this conjecture is at least valid up to
fourth order in the deformation.\par
%
%---------------------------------------------------------------------------
%
For such purpose, we have first to determine the commutators of $\bi{L}^{2k}$
and $\bi{A}^2$ with~$A$. We state the results in the form of two lemmas.\par
\bigskip
\noindent{\it Lemma 1.} For any $k\in\N^+$, the generators of the deformed
algebra~$g_d$, defined in~(3.1), satisfy the relation
$$\bigl[A, \bi{L}^{2k}\bigr]^1 = \sum_{i=0}^{k-1} x^{(k)}_i \bi{L}^{2i} A
    + \sum_{i=0}^{k-1} y^{(k)}_i \bi{L}^{2i} [L \times A]^1 + \sum_{i=0}
    ^{k-2} z^{(k)}_i \bi{L}^{2i} \left[[L \times L]^2 \times A\right]^1
    \eqno(3.11)$$
where $x^{(k)}_i$, $y^{(k)}_i$, $i=0$, 1, $\ldots$,~$k-1$, and $z^{(k)}_i$,
$i=0$, 1, $\ldots$,~$k-2$, are some real constants fulfilling the recursion
relations
$$\openup 2mm
  \eqalignno{x^{(k)}_i &= 2 x^{(k-1)}_i + x^{(k-1)}_{i-1} + \case{2}{3}
       \sqrt{2} y^{(k-1)}_{i-1} + 2 \delta_{i,k-1} \qquad i=0,1,\ldots,k-1
       &(3.12a) \cr
    y^{(k)}_i &= 2 \sqrt{2} x^{(k-1)}_i + y^{(k-1)}_i - \sqrt{\case{3}{10}}
       z^{(k-1)}_i + y^{(k-1)}_{i-1} + 2 \sqrt{\case{2}{15}} z^{(k-1)}_{i-1}
       \cr
    &\phantom{=} + 2 \sqrt{2} \delta_{i,k-1} \qquad i=0,1,\ldots,k-1 &(3.12b)
       \cr
    z^{(k)}_i &= \sqrt{\case{10}{3}} y^{(k-1)}_i - z^{(k-1)}_i +
       z^{(k-1)}_{i-1} \qquad i=0,1,\ldots,k-2 &(3.12c) \cr}$$
and the conditions $x^{(1)}_0 = 2$, $y^{(1)}_0 = 2 \sqrt{2}$.\par
\bigskip
\noindent{\it Lemma 2.} The generators of the deformed algebra~$g_d$, defined
in~(3.1), satisfy the relation
$$\bigl[A, \bi{A}^2\bigr]^1 = - \sum_{i=0}^{\infty} u_i \bi{L}^{2i} A
    - \sum_{i=0}^{\infty} v_i \bi{L}^{2i} [L \times A]^1 - \sum_{i=0}
    ^{\infty} w_i \bi{L}^{2i} \left[[L \times L]^2 \times A\right]^1
    \eqno(3.13)$$
where $u_i$, $v_i$, $w_i$ are some real constants defined by the formal
series
$$\eqalignno{u_i &= \sum_{k=i}^{\infty} a_k \left(2 x^{(k)}_i + \case{1}{3}
       \sqrt{2} y^{(k)}_{i-1} + 2 \delta_{k,i}\right) &(3.14a) \cr
    v_i &= \sum_{k=i}^{\infty} a_k \left(\sqrt{2} x^{(k)}_i + \case{3}{2}
       y^{(k)}_i - \case{1}{2} \sqrt{\case{3}{10}} z^{(k)}_i
       + \sqrt{\case{2}{15}} z^{(k)}_{i-1} + 2 \sqrt{2} \delta_{k,i}\right)
       &(3.14b) \cr
    w_i &= \sum_{k=i+1}^{\infty} a_k \left(\sqrt{\case{5}{6}} y^{(k)}_i
       + \case{1}{2} z^{(k)}_i\right) &(3.14c) \cr}$$
in terms of the solution of eq.~(3.12).\par
\bigskip
\noindent The proofs of lemmas~1 and~2 are sketched in appendix~2, and the
solution obtained for $x^{(k)}_i$, $y^{(k)}_i$, $z^{(k)}_i$,~$k\le 5$, by
solving eq.~(3.12), is listed in table~1.\par
%
%--------------------------------------------------------------------------
%
Finally, by combining lemmas~1 and~2, the following result can be easily
derived (for details see appendix~2):\par
\bigskip
\noindent{\it Proposition 3.} Up to fourth order in the deformation, the
operator~$C_{1d}$ defined in eq.~(3.10), where
$$\eqalignno{b_1 &= a_0 + a_1 \qquad b_2 = \case{1}{2} a_1 + \case{4}{3} a_2
       - \case{1}{3} a_3 + \case{8}{15} a_4 \qquad b_3 = \case{1}{3} a_3
       + \case{5}{3} a_3 - \case{16}{15} a_4 \cr
    b_4 &= \case{1}{4} a_3 + 2 a_4 \qquad b_5 = \case{1}{5} a_4 \qquad b_6 =
       b_7 = \cdots = 0 &(3.15)}$$
is a Casimir operator of the algebra~$g_d$, defined in~(3.1).\par
%
%xxxxxxxxxxxxxxxxxxxxxxxxxxxxxxxxxxxxxxxxxxxxxxxxxxxxxxxxxxxxxxxxxxxxxxxxxx
%
\subsection{Unitary irreducible representations}
In the present subsection, we will study the deformations of some unitary
irreducible representations (unirreps) of the Lie algebras $g = \hbox{\rm
so(4)}$, so(3,1), and~e(3) when the latter are replaced by the corresponding
deformed algebras~$g_d$.\par
%
%------------------------------------------------------------------------
%
The unirreps considered are those whose representation space $\cal R$
contains each of the representation spaces ${\cal R}^l$, $l=0$, $1\over2$, 1,
$3\over2$,~$\ldots$, of so(3) at most once (Biedenharn~1961, Naimark~1964,
B\H ohm~1979). Such unirreps can be characterized by
{\parindent= 1cm
\item{(i)} $[p,q]$ where $p\ge |q|$, and $p$, $|q|\in\N$ or $p$, $|q|\in
{1\over2}\N$, in the so(4) case,
\item{(ii)} $(l_0,c)$ where either $l_0 \in \left\{0,{1\over2},1,{3\over2},
\ldots\right\}$ and $c\in\R$, or $l_0=0$ and $c=i\nu$, $\nu\in\R$, in the
so(3,1) case,
\item{(iii)} $(l_0,\epsilon)$ where $l_0 \in \left\{0,{1\over2},1,{3\over2},
\ldots\right\}$ and $\epsilon\in\R$, in the e(3) case,\par}
\noindent where the latter are obtained from those of so(3,1) by an In\" on\"
u-Wigner contraction. In cases~(i) and~(ii) (or~(iii)), the decomposition of
their representation space is given by
$${\cal R} = \sum_{l=l_0,l_0+1}^{l_1} \oplus {\cal R}^l \qquad \hbox{\rm and}
    \qquad {\cal R} = \sum_{l=l_0,l_0+1}^{\infty} \oplus {\cal R}^l
    \eqno(3.16a,b)$$
respectively, where in~(3.16a), the minimum (resp.~maximum) $l$ value is
defined by $l_0 = |q|$ (resp.~$l_1 = p$).\par
%
%---------------------------------------------------------------------------
%
The reduced matrix elements of the vector operator~$A$ and the eigenvalues of
the Casimir operators can be written as\footnote\dag{The phase convention
adopted in the present paper is that of Biedenharn~(1961), which differs from
that of Naimark~(1964) and B\H ohm~(1979).}
$$\eqalignno{\langle[p,q]l \|A\| [p,q]l\rangle &= {q(p+1) \over
       [l(l+1)]^{1/2}} \cr
    \langle[p,q]l-1 \|A\| [p,q]l\rangle &= - \left[{(l-q) (l+q) (p+1-l)
       (p+1+l) \over l(2l-1)}\right]^{1/2} &(3.17) \cr}$$
and
$$\langle C_1 \rangle = p(p+2) + q^2 \qquad \langle C_2 \rangle = q(p+1)
    \eqno(3.18)$$
for so(4),
$$\eqalignno{\langle(l_0,c)l \|A\| (l_0,c)l\rangle &= - {l_0 c \over
       [l(l+1)]^{1/2}} \cr
    \langle(l_0,c)l-1 \|A\| (l_0,c)l\rangle &= - \left[{(l-l_0) (l+l_0)
       (c^2+l^2) \over l(2l-1)}\right]^{1/2} &(3.19) \cr}$$
and
$$\langle C_1 \rangle = c^2 - l_0^2 + 1 \qquad \langle C_2 \rangle = - l_0 c
    \eqno(3.20)$$
for so(3,1),
$$\eqalignno{\langle(l_0,\epsilon)l \|A\| (l_0,\epsilon)l\rangle &= - {l_0
       \epsilon \over [l(l+1)]^{1/2}} \cr
    \langle(l_0,\epsilon)l-1 \|A\| (l_0,\epsilon)l\rangle &= - |\epsilon|
       \left[{(l-l_0) (l+l_0) \over l(2l-1)}\right]^{1/2} &(3.21) \cr}$$
and
$$\langle C_1 \rangle = \epsilon^2 \qquad \langle C_2 \rangle = - l_0
    \epsilon \eqno(3.22)$$
for e(3). In all cases, the remaining reduced matrix elements of $A$ can be
obtained from the relation
$$\langle l+1 \|A\| l\rangle = - \left({2l+1 \over 2l+3}\right)^{1/2}
    \langle l \|A\| l+1\rangle \eqno(3.23)$$
valid for any vector operator.\par
%
%---------------------------------------------------------------------------
%
By introducing the function $G(l^2,l_0^2)$, defined by
$$\eqalignno{G(l^2,l_0^2) &= {1 \over l^2 - l_0^2} \sum_{j=l_0}^{l-1} (2j+1)
       g\bigl(j(j+1)\bigr) \cr
    &= a_0 + \case{1}{2} a_1 (l^2 + l_0^2 - 1) + \case{1}{3} a_2
       \left[(l^2-1)^2 + l_0^2 (l^2-1) + l_0^2 (l_0^2-1)\right] \cr
    &\phantom{=} + \cdots &(3.24) \cr}$$
for $l>l_0$, we easily find the following results:\par
\bigskip
\noindent{\it Proposition 4.} Provided the deformation parameters $a_k$,
$k>0$, are chosen in such a way that all quantities under square roots remain
nonnegative, the unirreps~$[p,q]$, $(l_0,c)$, and~$(l_0,\epsilon)$ of so(4),
so(3,1), and~e(3), respectively, can be deformed into unirreps of the
corresponding deformed algebras~$g_d$. The reduced matrix elements of~$A$ and
the eigenvalues of the Casimir operators become
$$\eqalign{\langle[p,q]l \|A\| [p,q]l\rangle &= q(p+1) \left[{ G\left(
       (p+1)^2,q^2\right) \over l(l+1)}\right]^{1/2} \cr
    \langle[p,q]l-1 \|A\| [p,q]l\rangle &= - \left\{{(l-q) (l+q) \left[
       (p+1)^2 G\left((p+1)^2,q^2\right) - l^2 G\left(l^2,q^2\right)\right]
       \over l(2l-1)} \right\}^{1/2} \cr} \eqno (3.25)$$
and
$$\eqalignno{\langle C_{1d} \rangle &= h\bigl(|q|(|q|+1)\bigr) - (|q|+1)
       G\left((|q|+1)^2,q^2\right) + (p+1)^2 G\left((p+1)^2,q^2\right) \cr
    &= \langle C_1 \rangle + \case{1}{2} a_1 \Bigl[\langle C_1 \rangle
       \left(\langle C_1 \rangle + 1\right) - \langle C_2 \rangle^2 \Bigr]
       + \case{1}{3} a_2 \langle C_1 \rangle \Bigl[\langle C_1 \rangle
       \left(\langle C_1 \rangle + 1\right) - 2 \langle C_2 \rangle^2 \Bigr]
       \cr
    &\phantom{=} + \cdots \cr
    \langle C_{2d} \rangle &= q(p+1) \Bigl[G\bigl((p+1)^2,q^2\bigr)\Bigr]
       ^{1/2} \cr
    &= \langle C_2 \rangle \Bigl[1 + \case{1}{2} a_1 \langle C_1 \rangle
       + \case{1}{3} a_2 \left(\langle C_1 \rangle^2 - \langle C_2 \rangle^2
       \right) + \cdots\Bigr]^{1/2} &(3.26) \cr}$$
in the deformed so(4) case,
$$\eqalignno{\langle(l_0,c)l \|A\| (l_0,c)l\rangle &= - {l_0 c \over
       [l(l+1)]^{1/2}} \cr
    \langle(l_0,c)l-1 \|A\| (l_0,c)l\rangle &= - \left\{{(l-l_0) (l+l_0)
       \left[c^2 - l^2 G(l^2,l_0^2)\right] \over l(2l-1)}\right\}^{1/2}
       &(3.27) \cr}$$
and
$$\eqalignno{\langle C_{1d} \rangle &= h\bigl(l_0(l_0+1)\bigr) - (l_0+1)
       G\left((l_0+1)^2,l_0^2\right) + c^2 \cr
    &= \langle C_1 \rangle + \case{1}{2} a_1 l_0^2 \left(l_0^2-1\right)
       + \case{1}{3} a_2  l_0^2 \left(l_0^2-1\right)^2 + \cdots \cr
    \langle C_{2d} \rangle &=  \langle C_2 \rangle = - l_0 c &(3.28) \cr}$$
in the deformed so(3,1) case, while for deformed e(3) they can be obtained by
substituting $\epsilon$ for $c$ in eqs.~(3.27) and~(3.28).\par
\bigskip
\noindent{\it Proof.} Whenever all $a_k$'s, for which $k>0$, go to zero,
equations~(3.25) to~(3.28), and their counterparts for e(3) go over into the
undeformed results contained in eqs.~(3.17) to~(3.22), respectively. On the
other hand, for arbitrary values of the $a_k$'s satisfying the hypothesis,
the validity of the equations can be checked by direct substitution into the
commutation relation~(3.1c), and the definitions of $C_{1d}$
and~$C_{2d}$.\cqfd\par
\bigskip
\noindent{\it Remark.} For some choices of the deformation parameters, it may
happen that the unirreps of $g_d$ considered in proposition~4 do not exhaust
the class of unirreps whose representation space contains each of the
representation spaces of so(3) at most once. The existence of ``extra''
representations, which have no counterpart for the undeformed algebra, has
already been noted in the deformed su(2) case (Ro\v cek~1991).\par
%
%==========================================================================
%
\section{The quadrupole case}
In the quadrupole case, the deformed algebra~$g_d$ is defined by the
commutation relations
$$\eqalignno{[L,L]^1 &= -\sqrt{2} L \qquad [L,Q]^{\Lambda} = - \sqrt{6}\,
        \delta_{\Lambda,2} Q &(4.1a,b) \cr
   [Q,Q]^1 &= f^1(L) = 3 \sqrt{10}\, g_1(\bi{L}^2) L = 3 \sqrt{10} \biggl(
        \sum_{k=0}^{\infty} a^{(1)}_k \bi{L}^{2k}\biggr) L &(4.1c) \cr
   [Q,Q]^3 &= f^3(L) = g_3(\bi{L}^2) \left[[L \times L]^2 \times L
        \right]^3 \cr
   &=  \biggl( \sum_{k=0}^{\infty} a^{(3)}_k \bi{L}^{2k}\biggr) \left[[L
        \times L]^2 \times L\right]^3 &(4.1d) \cr}$$
where the normalization constant~$\gamma_1$ has been set equal to~$3
\sqrt{10}$. The algebra is associative provided $Q$ satisfies eqs.~(2.7a,b).
\par
%
%--------------------------------------------------------------------------
%
In the undeformed case, where the algebra~$g_d$ reduces to an ordinary Lie
algebra~$g$, one has $a^{(1)}_k=0$, $k=1$, 2,~$\ldots$, and $a^{(3)}_k=0$,
$k=0$, 1,~$\ldots$. According to whether $a^{(1)}_0=+1$, $-1$, or~0, $g$
is the special unitary algebra su(3) (Elliott~1958a,b), the special linear
algebra sl(3,\R) (Weaver and Biedenharn~1972), or the semidirect sum of an
abelian algebra t(5) and so(3) (Ui~1970, Weaver \etal 1973).\par
%
%--------------------------------------------------------------------------
%
If we restrict ourselves to a first-order deformation, equations~(4.1c)
and~(4.1d) become
$$\eqalignno{[Q,Q]^1 &= 3 \sqrt{10}\, \epsilon L + \alpha \bi{L}^2 L &(4.2a)
\cr
       [Q,Q]^3 &= \beta \left[[L \times L]^2 \times L\right]^3 &(4.2b) \cr}$$
where $\alpha$, $\beta$, and $\epsilon$ are defined by $\alpha = 3 \sqrt{10}
a^{(1)}_1$, $\beta = a^{(3)}_0$, and $\epsilon = a^{(1)}_0 = +1$, $-1$,
or~0. By inserting eqs.~(4.2a,b) into eqs.~(2.7a,b) and taking eq.~(4.1b)
into account, we obtain the two associativity conditions
$$\eqalignno{2 \sqrt{2}\, \alpha \left[Q, \bi{L}^2 L\right]^1 + \sqrt{7}\,
            \beta \left[Q, \left[[L \times L]^2 \times L\right]^3 \right]^1
            &= 0 &(4.3a) \cr
       \alpha \left[Q, \bi{L}^2 L\right]^3 - 2 \beta \left[Q, \left[[L
            \times L]^2 \times L\right]^3 \right]^3 &= 0. &(4.3b) \cr}$$\par
%
%--------------------------------------------------------------------------
%
Straightforward tensor algebra leads to the following results:
$$\eqalignno{\left[Q, \bi{L}^2 L\right]^{\Lambda} &= 2 \sqrt{3} \Bigl\{
            \bigl[\sqrt{3} - U(11\Lambda2;12) \bigr] [L \times Q]^{\Lambda}
            \cr
       &\phantom{=} + \sqrt{2}\, U(11\Lambda2;22) \left[[L \times L]^2
            \times Q \right]^{\Lambda} \Bigr\} \cr
      \left[Q, \left[[L \times L]^2 \times L\right]^3 \right]^{\Lambda} &=
            \sqrt{3} \Bigl\{ \bigl[\sqrt{7}\, U(22\Lambda1;23) \cr
       &\phantom{=} - 2 \sqrt{2\Lambda(\Lambda+1)}\, U(22\Lambda1;\Lambda3)
           U(21\Lambda1;22) \bigr] [L \times Q]^{\Lambda} \cr
       &\phantom{=} + 3 \sqrt{2}\, U(21\Lambda2;23) \left[[L \times L]^2
           \times Q\right]^{\Lambda} \Bigr\} &(4.4)\cr}$$
valid for $\Lambda=1$ or~3. By replacing Racah coefficients by their
numerical values in eq.~(4.4), conditions~(4.3a,b) can be rewritten as
$$\eqalignno{\bigl(6 \sqrt{10}\, \alpha - 7 \beta\bigr) [L \times Q]^1 +
       \sqrt{6} \bigl(2 \sqrt{10}\, \alpha + 7 \beta\bigr) \left[[L \times
       L]^2 \times Q\right]^1 &= 0 &(4.5a) \cr
    2 \bigl(\sqrt{10}\, \alpha + 3 \beta\bigr) [L \times Q]^3 +
       \bigl(\sqrt{10}\, \alpha - 9 \beta\bigr) \left[[L \times L]^2 \times
       Q \right]^3 &= 0. &(4.5b) \cr}$$
Since they cannot be satisfied for any choice of the deformation parameters
$\alpha$,~$\beta$, we conclude that for a first-order deformation, the
algebra~$g_d$, defined in~(4.1), is not associative contrary to what happens
for the algebra~(3.1) in the vector case. We shall therefore pursue the
analysis of the quadrupole case no further.\par
%
%==========================================================================
%
\section{Conclusion}
In the present paper, we established that there exist deformations of the
Lie algebras so(4), so(3,1), and~e(3) that leave their so(3) subalgebra
undeformed and preserve their coset structure. We proved that the Casimir
operators of these Lie algebras can be deformed so as to provide
corresponding operators for the deformed algebras. Moreover, we constructed
those unirreps of the latter that go over for vanishing deformation into the
unirreps of the former belonging to an important class of representations.
\par
%
%--------------------------------------------------------------------------
%
In contrast, we showed that a similar deformation of the Lie algebras su(3),
sl(3,\R), and of the semidirect sum of t(5) and so(3) is not possible
because the associativity conditions are violated in first order in the
deformation.\par
%
%--------------------------------------------------------------------------
%
It should be stressed that the deformations of so(4), so(3,1),
and e(3) studied here differ from the standard $q$-algebras so$_q$(4),
so$_q$(3,1), e$_q$(3) (Drinfeld~1986, Jimbo~1985, Celeghini \etal 1991,
Chakrabarti~1993), as well as from an alternative deformation of the
orthogonal and pseudo-orthogonal Lie algebras proposed by Gavrilik and
Klimyk~(1991) (see also Gavrilik~1993). In both these approaches, the so(3)
subalgebra is indeed deformed contrary to what happens in the present one.
Since rotational invariance and angular momentum conservation are important
properties of many physical systems, one may hope that the deformed algebras
introduced in this paper will prove more relevant to applications than
those previously considered.\par
%
%---------------------------------------------------------------------------
%
Some problems wherein the coset structure so(4)/so(3) is important can be
found in standard quantum mechanics (e.g., that of a particle in a Coulomb
potential (Biedenharn~1961)), as well as in parasupersymmetric quantum
mechanics with three parasupercharges (Debergh and Nikitin~1995).
Deformations of so(4) preserving that coset structure may therefore be
expected to play a role in similar contexts. It is already known that the
first-order deformation of so(4) is the symmetry algebra of a particle in a
Coulomb potential when the space has a constant curvature (Higgs~1979,
Leemon~1979, Granovskii \etal 1992, de Vos and van Driel~1993). All the
general results derived in the present paper therefore apply to such a
problem. At a more phenomenological level, deviations from hydrogenic
spectra that are found for many-electron atoms or excitons in semiconductors
might be accounted for by some deformations of so(4). Similarly,
deformations of parasupersymmetric quantum mechanics with three
parasupercharges might lead to some parasupersymmetric Hamiltonians with new
and nontrivial properties, as happens in the case of two parasupercharges
(Debergh \etal 1995). We hope to come back to some of these problems in
forthcoming publications.\par
\vfill\eject
%
%=========================================================================
%
\appendix{Definition and properties of coupled commutators}
The purpose of this appendix is to review the definition and some useful
properties of coupled commutators.\par
%
%-------------------------------------------------------------------------
%
The coupled commutator of two so(3) irreducible tensors $T^{\lambda_1}$ and
$U^{\lambda_2}$, of ranks $\lambda_1$ and $\lambda_2$ respectively, is
defined by
$$\left[T^{\lambda_1}, U^{\lambda_2} \right]^{\Lambda}_M = \sum_{\mu_1
       \mu_2} \langle \lambda_1 \mu_1, \lambda_2 \mu_2 | \Lambda M \rangle
       \left[T^{\lambda_1}_{\mu_1}, U^{\lambda_2}_{\mu_2} \right] \eqno({\rm
       A}1.1)$$
in terms of an ordinary commutator $[\, , \,]$ and of an su(2) Wigner
coefficient $\langle\, , \,|\, \rangle$. By using a symmetry property of the
latter (Rose~1957), equation~(A1.1) can be alternatively written as
$$\eqalignno{\left[T^{\lambda_1}, U^{\lambda_2} \right]^{\Lambda}_M &=
            \left[T^{\lambda_1} \times U^{\lambda_2} \right]^{\Lambda}_M -
            (-1)^{\lambda_1 + \lambda_2 - \Lambda} \left[U^{\lambda_2} \times
            T^{\lambda_1} \right]^{\Lambda}_M \cr
       &= - (-1)^{\lambda_1 + \lambda_2 - \Lambda} \left[U^{\lambda_2},
            T^{\lambda_1} \right]^{\Lambda}_M &({\rm A}1.2) \cr}$$
where
$$\left[T^{\lambda_1} \times U^{\lambda_2} \right]^{\Lambda}_M =
       \sum_{\mu_1 \mu_2} \langle \lambda_1 \mu_1, \lambda_2 \mu_2 |
       \Lambda M \rangle \, T^{\lambda_1}_{\mu_1} U^{\lambda_2}_{\mu_2}.
       \eqno({\rm A}1.3)$$\par
%
%--------------------------------------------------------------------------
%
For three irreducible tensors $T^{\lambda_1}$, $U^{\lambda_2}$,
$V^{\lambda_3}$, the well-known relation $[A, BC] = [A, B] C + B [A, C]$
becomes in coupled form
$$\eqalignno{\Bigl[T^{\lambda_1}, \left[U^{\lambda_2} \times V^{\lambda_3}
            \right]^{\Lambda_{23}} \Bigr]^{\Lambda}_M &= \sum_{\Lambda_{12}}
            U(\lambda_1 \lambda_2 \Lambda \lambda_3; \Lambda_{12}
            \Lambda_{23}) \left[ \left[T^{\lambda_1}, U^{\lambda_2}\right]
            ^{\Lambda_{12}}  \times V^{\lambda_3} \right]^{\Lambda}_M \cr
       &\phantom{=} + \sum_{\Lambda_{13}} (-1)^{\lambda_3 + \Lambda -
            \Lambda_{13} - \Lambda_{23}} U(\lambda_1 \lambda_3 \Lambda
            \lambda_2;\Lambda_{13} \Lambda_{23}) \cr
       &\phantom{=} \times \left[U^{\lambda_2} \times \left[
           T^{\lambda_1}, V^{\lambda_3} \right]^{\Lambda_{13}}
           \right]^{\Lambda}_M &({\rm A}1.4)}$$
where $U(\lambda_1 \lambda_2 \Lambda \lambda_3; \Lambda_{12}
\Lambda_{23})$ denotes a Racah coefficient in unitary form (Rose~1957). In
the same way, the Jacobi identity $\bigl[A, [B, C]\bigr] + \bigl[B, [C, A]
\bigr] + \bigl[C, [A, B]\bigr] = 0$ can be rewritten as
$$\deqn{\left[T^{\lambda_1}, \left[U^{\lambda_2}, V^{\lambda_3} \right]
           ^{\Lambda_{23}} \right]^{\Lambda}_M + \sum_{\Lambda_{31}}
           (-1)^{\lambda_1 + \Lambda_{23} - \Lambda} U(\lambda_2 \lambda_3
           \Lambda \lambda_1; \Lambda_{23} \Lambda_{31}) \left[U^{\lambda_2},
           \left[V^{\lambda_3}, T^{\lambda_1} \right]^{\Lambda_{31}}
           \right]^{\Lambda}_M \cr
      \ind + \sum_{\Lambda_{12}}
           (-1)^{\lambda_3 + \Lambda_{12} - \Lambda} U(\lambda_1 \lambda_2
           \Lambda \lambda_3; \Lambda_{12} \Lambda_{23}) \left[V^{\lambda_3},
           \left[T^{\lambda_1}, U^{\lambda_2} \right]^{\Lambda_{12}}
           \right]^{\Lambda}_M \cr
      \ind \eql 0. \eq({\rm A}1.5)}$$\par
\vfill\eject
%
%===========================================================================
%
\appendix{Proofs of lemmas~1 and~2 and of proposition~3}
To prove eqs.~(3.11) and~(3.12), we proceed by induction over $k$. For the
lowest value $k=1$, equation~(3.11) reduces to (3.4). For $k>1$, we start
from the identities
$$\bigl[A, \bi{L}^{2k}\bigr]^1 = \bi{L}^{2k-2} \left[A, \bi{L}^2\right]^1 +
    \bi{L}^2 \bigl[A, \bi{L}^{2k-2}\bigr]^1 + \Bigl[ \bigl[A,
    \bi{L}^{2k-2}\bigr]^1, \bi{L}^2 \Bigr]^1 \eqno({\rm A}2.1)$$
and
$$\Bigl[ \bigl[A, \bi{L}^{2k-2}\bigr]^1, \bi{L}^2 \Bigr]^1 = 2 \biggl\{
      \bigl[A, \bi{L}^{2k-2}\bigr]^1 + \sqrt{2} \Bigl[L \times \bigl[A,
      \bi{L}^{2k-2}\bigr]^1 \Bigr]^1 \biggr\} \eqno({\rm A}2.2)$$
resulting from eq.~(A1.4) and standard tensor algebra. Then assuming
eq.~(3.11) to be valid when $k$ is replaced by $k-1$, we use it to compute
the right-hand sides of eqs.~(A2.1) and~(A2.2). The result contains two
tensor products that are not in standard form, as those appearing on the
right-hand side of~(3.11), but can be rewritten in such a form by using the
following identities:
$$\left[L \times [L \times A]^1 \right]^1 = \case{1}{3} \bi{L}^2 A -
       \case{1}{2 \sqrt{2}} [L \times A]^1 + \case{1}{2} \sqrt{\case{5}{3}}
       \left[[L \times L]^2 \times A \right]^1 \eqno({\rm A}2.3)$$
$$\left[L \times \left[ [L \times L]^2 \times A \right]^1 \right]^1 = -
       \case{1}{4\sqrt{15}} \left(3 - 4 \bi{L}^2\right) [L \times A]^1 -
       \case{3}{2\sqrt{2}} \left[[L \times L]^2 \times A \right]^1.
       \eqno({\rm A}2.4)$$
After some straightforward calculations, equation~(3.11) is obtained provided
$x^{(k)}_i$, $y^{(k)}_i$, and $z^{(k)}_i$ satisfy eq.~(3.12).\par
%
%--------------------------------------------------------------------------
%
The proof of eqs.~(3.13) and~(3.14) is based upon the identities
$$\left[A, \bi{A}^2\right] = - 2 \sqrt{2} \sum_{k=0}^{\infty} a_k \bi{L}^{2k}
    [L \times A]^1 + \sqrt{2} \sum_{k=0}^{\infty} a_k \bigl[A, \bi{L}^{2k} L
    \bigr]^1 \eqno({\rm A}2.5)$$
and
$$\bigl[A, \bi{L}^{2k} L\bigr]^1 = - \sqrt{2} \bigl[A, \bi{L}^{2k}\bigr]^1
    - \Bigl[L \times \bigl[A, \bi{L}^{2k}\bigr]^1 \Bigr]^1 - \sqrt{2}
    \bi{L}^{2k} A \eqno({\rm A}2.6)$$
resulting from eqs.~(A1.4) and~(3.1), as well as standard tensor algebra.
Taking eqs.~(3.11), (A2.3), and (A2.4) into account directly leads to the
searched for results.\par
%
%--------------------------------------------------------------------------
%
Finally, from eqs.~(3.10), (3.11), and (3.13), it follows that for arbitrary
constants~$a_k$, the condition $[A, C_{1d}]^1 = 0$ is fulfilled provided
$$\sum_{k=i+1}^{\infty} b_k x^{(k)}_i - u_i = \sum_{k=i+1}^{\infty} b_k
    y^{(k)}_i - v_i = \sum_{k=i+2}^{\infty} b_k z^{(k)}_i - w_i = 0 \qquad
    i=0,1,2,\ldots. \eqno({\rm A}2.7)$$
By successively using eqs.~(3.14) and~(3.12), these conditions can be
rewritten as
$$\eqalignno{2 \sum_{k=i+1}^{\infty} b'_k x^{(k)}_i &=
       \sum_{k=i}^{\infty} a_k \left(2 x^{(k)}_i - x^{(k)}_{i-1} +
       2\delta_{k,i} \right) \qquad i=0,1,2,\ldots &({\rm A}2.8a) \cr
    2 \sum_{k=i+1}^{\infty} b'_k y^{(k)}_i &=
       \sum_{k=i}^{\infty} a_k \left(2 y^{(k)}_i - y^{(k)}_{i-1} + 2 \sqrt{2}
       \delta_{k,i}\right) \qquad i=0,1,2,\ldots &({\rm A}2.8b) \cr
    2 \sum_{k=i+2}^{\infty} b'_k z^{(k)}_i &=
       \sum_{k=i+1}^{\infty} a_k \left(2 z^{(k)}_i - z^{(k)}_{i-1} \right)
       \qquad i=0,1,2,\ldots &({\rm A}2.8c) \cr}$$
where $b'_k \equiv b_k - {1\over2} a_{k-1}$.\par
%
%---------------------------------------------------------------------------
%
If we restrict ourselves to a $K$th-order deformation and assume that $b_k=0$
for $k>K+1$, conditions~(A2.8a) and~(A2.8b) reduce to two systems of $K+1$
equations (corresponding to $i=0$, 1, $\ldots$,~$K$) in $K+1$ unknowns
$b'_k$, $k=1$, 2, $\ldots$,~$K+1$, while condition~(A2.8c) leads to a system
of $K$ equations (corresponding to $i=0$, 1, $\ldots$,~$K-1$) in $K$
unknowns $b'_k$, $k=2$, 3, $\ldots$,~$K+1$. It only remains to solve (A2.8a)
and to check that its solution also satisfies the two remaining systems of
equations. This calculation was carried out for $K=4$, and the results are
contained in eq.~(3.15).\par
%
%==========================================================================
%
\references
\refbk{Beckers J, Debergh N and Quesne C 1995 Parasupersymmetric quantum
mechanics with generalized deformed parafermions}{Preprint}{ULB-ULG}
\refjl{Biedenharn L C 1961}{\JMP}{2}{433}
\refbk{B\H ohm A 1979}{Quantum Mechanics}{(Berlin: Springer) p 146}
\refbk{Bonatsos D, Daskaloyannis C, Kolokotronis P and Lenis D 1994
Nonlinear extension of the u(2) algebra as the symmetry algebra of the
planar anisotropic quantum harmonic oscillator with rational ratio of
frequencies and ``pancake'' nuclei}{Preprint}{NUCL-TH/9412003}
\refjl{Celeghini E, Giachetti R, Sorace E and Tarlini M
1991}{\JMP}{32}{1159} %
\refjl{Chakrabarti A 1993}{\JMP}{34}{1964}
\refbk{Debergh N and Nikitin A G 1995 Parasupersymmetric quantum mechanics
with an arbitrary number of parasupercharges and orthogonal Lie
algebras}{Helv. Phys. Acta}{in press}
\refjl{de Boer J and Tjin T 1993}{Commun. Math. Phys.}{158}{485}
\refbk{de Vos K and van Driel P 1993 The Kazhdan-Lusztig conjecture for
finite W-algebras}{Preprint}{hep-th/9312016}
\refbk{Drinfeld V G 1986}{Proc. Int. Congr. of Mathematicians (Berkeley,
CA)}{ed.\ A M Gleason (Providence, RI: AMS) p 798}
\refjl{Elliott J P 1958a}{Proc. R. Soc.{\rm\ A}}{245}{128}
\refjl{\dash 1958b}{Proc. R. Soc.{\rm\ A}}{245}{562}
\refjl{Fairlie D B and Nuyts J 1994}{\JMP}{35}{3794}
\refjl{Fairlie D B and Zachos C K 1991}{\PL}{256B}{43}
\refjl{Gaskell R, Peccia A and Sharp R T 1978}{\JMP}{19}{727}
\refjl{Gavrilik A M 1993}{Theor. Math. Phys.}{95}{546}
\refjl{Gavrilik A M and Klimyk A U 1991}{Lett. Math. Phys.}{21}{215}
\refjl{Granovskii Ya I, Zhedanov A S and Lutsenko I M 1992}{Theor. Math.
Phys.}{91}{604}
\refjl{Higgs P W 1979}{\JPA}{12}{309}
\refjl{Jimbo M 1985}{Lett. Math. Phys.}{10}{63}
\refjl{Karassiov V P 1994}{\JPA}{27}{153}
\refjl{Karassiov V P and Klimov A B 1994}{\PL}{191A}{117}
\refjl{Leemon H I 1979}{\JPA}{12}{489}
\refbk{Naimark M A 1964}{Linear Representations of the Lorentz Group}{(New
York: Pergamon)}
\refjl{Polychronakos A P 1990}{Mod. Phys. Lett.{\rm\ A}}{5}{2325}
\refjl{Quesne C 1994}{\PL}{193A}{245}
\refjl{Ro\v cek M 1991}{\PL}{255B}{554}
\refbk{Rose M E 1957}{Elementary Theory of Angular Momentum}{(New York:
Wiley)}
\refjl{Schoutens K, Sevrin A and van Nieuwenhuizen P 1989}{Commun. Math.
Phys.}{124}{87}
\refjl{Tjin T 1992}{\PL}{292B}{60}
\refjl{Ui H 1970}{Prog. Theor. Phys.}{44}{153}
\refjl{Weaver L and Biedenharn L C 1972}{\NP{\rm\ A}}{185}{1}
\refjl{Weaver L, Biedenharn L C and Cusson R Y 1973}{\APNY}{77}{250}
\refjl{Zamolodchikov A B 1986}{Theor. Math. Phys.}{65}{1205}
%
%=========================================================================
%
\tables

\tabcaption{Solution of the recursion relations~(3.12) up to $k=5$.}
\bigskip
\catcode`\*=\active
\def*{\phantom{-}}
\vbox{\openup 2mm
\halign{{$\displaystyle#$}\hfill&&\qquad{$\displaystyle#$}\hfill\cr
        \noalign{\boldrule{22pc}\vskip -2mm}
        \omit & i=0& i=1& i=2& i=3 &i=4\cr
        \noalign{\vskip -2mm\medrule{22pc}\vskip -2mm}
        x^{(2)}_i& 4& \case{20}{3}& \bar& \bar& \bar\cr
        y^{(2)}_i& 6\sqrt{2}& 4\sqrt{2}& \bar& \bar& \bar\cr
        z^{(2)}_i& 4\sqrt{\case{5}{3}}& \bar& \bar& \bar& \bar\cr
        x^{(3)}_i& 8& \case{76}{3}& 14& \bar& \bar\cr
        y^{(3)}_i& 12\sqrt{2}& 26\sqrt{2}& 6\sqrt{2}& \bar& \bar\cr
        z^{(3)}_i& 8\sqrt{\case{5}{3}}& 4\sqrt{15}& \bar& \bar& \bar\cr
        x^{(4)}_i& 16& \case{224}{3}& 88& 24& \bar\cr
        y^{(4)}_i& 24\sqrt{2}& 88\sqrt{2}& 68\sqrt{2}& 8\sqrt{2}& \bar\cr
        z^{(4)}_i& 16\sqrt{\case{5}{3}}& 16\sqrt{15}& 8\sqrt{15}& \bar& \bar
          \cr
        x^{(5)}_i& 32& \case{592}{3}& 368& \case{680}{3}& \case{110}{3}\cr
        y^{(5)}_i& 48\sqrt{2}& 248\sqrt{2}& 352\sqrt{2}& 140\sqrt{2}&
          10\sqrt{2}\cr
        z^{(5)}_i& 32\sqrt{\case{5}{3}}& 48\sqrt{15}&
          160\sqrt{\case{5}{3}}& 40\sqrt{\case{5}{3}}& \bar\cr
        \noalign{\vskip -2mm\boldrule{22pc}}
        }}
\vfill\eject

\bye